\documentclass[11pt]{article}
\usepackage{a4wide}

\usepackage[a4paper, margin=2.3cm]{geometry}
\usepackage{amsmath,amssymb,amsthm, bm}
\usepackage{color}
\usepackage{parskip}
\usepackage{hyperref}
\usepackage{parskip}
\usepackage{setspace}
\usepackage{graphicx}
\numberwithin{equation}{section}

\usepackage{setspace}
\usepackage{graphicx}
\usepackage{mathtools}

\theoremstyle{plain}
\newtheorem{theorem}{Theorem}[section]

\newtheorem{corollary}[theorem]{Corollary}

\theoremstyle{definition}

\theoremstyle{remark}
\newtheorem{remark}[theorem]{Remark}

\newtheorem{case[theorem]}{Case}

\usepackage{fixme}
\fxsetup{
    status=draft,
    author=,
    layout=margin,
    theme=color
}

\def\P{\mathcal{P}}
\def\L{\mathcal{L}}

\def\F{{\mathbb F}}
\def\R{{\mathbb R}}

\def\C{{\mathbb C}}

\begin{document}
\title{A Point-Conic Incidence Bound and Applications over $\F_p$}

\author{Ali Mohammadi \thanks{School of Mathematics, Institute for Research in Fundamental Sciences (IPM), Tehran, Iran. Email: \href{mailto:a.mohammadi@ipm.ir}{a.mohammadi@ipm.ir}}\and Thang Pham \thanks{Department of Mathematics, HUS, Vietnam National University. Email: \href{mailto:thangpham.math@vnu.edu.vn}{thangpham.math@vnu.edu.vn}} \and Audie Warren \thanks{Johann Radon Institute for Computational and Applied Mathematics (RICAM), Linz, Austria. Email: \href{mailto:audie.warren@oeaw.ac.at}{audie.warren@oeaw.ac.at}}}
\date{}
\maketitle 
\begin{abstract}
In this paper, we prove the first incidence bound for points and conics over prime fields. As applications, we prove new results on expansion of bivariate polynomial images and on certain variations of distinct distances problems. These include new lower bounds on the number of pinned algebraic distances as well as improvements of results of Koh and Sun (2014) and Shparlinski (2006) on the size of the distance set formed by two large subsets of finite dimensional vector spaces over finite fields. We also prove a variant of Beck's theorem for conics. 
\end{abstract}
\section{Introduction}
For an arbitrary field $\F$ and sets of points $\P$ and algebraic curves $\mathcal{C}$ in $\F^d$, we denote the number of incidences between $\P$ and $\mathcal{C}$ by $I(\P, \mathcal{C}) := |\{(\bm{p}, C)\in \P\times \mathcal{C}: \bm{p}\in C\}|$.
The celebrated Szemer\'{e}di--Trotter theorem \cite{sze-tt} states that for finite sets of points $\P$ and lines $\L$ over $\R^2$, one has
\begin{equation}\label{szt}
I(\P, \L)\ll |\P|^{\frac{2}{3}}|\L|^{\frac{2}{3}}+|\P|+|\L|.
\end{equation}
It is well-known that this bound is sharp up to a constant factor. Many applications of this result can be found in \cite{Zdvir, tao}, and references therein.

Let $\F_q$ denote a finite field of order $q$ and characteristic $p$. A finite field analogue of the Szemer\'{e}di--Trotter theorem was first studied by Bourgain, Katz and Tao \cite{bkt} in 2004. They proved that for any point set $\P$ and any line set $\L$ in $\mathbb{F}_p^2$ with $|\P|=|\L|=N=p^{\alpha}$, $0<\alpha<2$, we have \begin{equation}\label{eqszt1}I(\P, \L)\ll N^{\frac{3}{2}-\varepsilon},~ \mbox{where}~ \varepsilon=\varepsilon(\alpha)>0.\end{equation} The study of this type of incidence structure is not only interesting from a geometric perspective, but is also largely motivated by a wide range of applications in different areas, such as arithmetic combinatorics, number theory, restriction theory, and theoretical computer science. Since its appearance nearly two decades ago, only a few quantitative variants of estimate (\ref{eqszt1}) have been proved. In particular, when $N$ is large with respect to the order of the field, say $N\ge p$, Vinh \cite{Vinh} showed that 
\begin{equation}\label{eqszt2}\left\vert I(\P, \L)-\frac{|\P||\L|}{p}\right\vert\le p^{1/2}\sqrt{|\P||\L|}.\end{equation}
This was proved using techniques from graph theory. A very short and elementary proof can also be found in \cite{MP1}. This result also holds for the setting of arbitrary finite fields $\mathbb{F}_q$. When $|\P|=|\L|=p^{3/2}$, it follows from (\ref{eqszt2}) that $I(\P, \L)\sim N^{4/3}$, which matches the bound of (\ref{szt}).  On the other hand, it has been indicated in \cite{vinh-sharp} that the lower bound of (\ref{eqszt2}) is sharp in the sense that there are sets $\P$ and $\L$ with $|\P||\L|\gg p^3$ and there is no incidence between $\P$ and $\L$.  When $N$ is extremely small, say, $N\ll \log_2\log_6\log_{18}p$, Grosu \cite{grosu} proved that the point-line incidence structure in $\mathbb{F}_p^2$ is almost the same as that in $\mathbb{C}^2$. As a consequence, relying on the $\C^2$ analogue of estimate~\eqref{szt} due to T\'oth~\cite{Toth}, he obtained $I(\P, \L)\ll N^{4/3}$. Therefore, we are left with the situation that $\log_2\log_6\log_{18}p\ll N\ll p^{3/2}$. 

In the case $|\P|=|\L|=N\le p$, Helfgott and Rudnev \cite{HH} proved that $I(\P, \L)\ll N^{\frac{3}{2}-\frac{1}{10678}}$. Jones \cite{jone} showed that the saving of $\frac{1}{10678}$ can be improved to $\frac{1}{662}$. Both results in \cite{HH, jone} were proved by using sum-product type energy inequalities in the field $\mathbb{F}_p$. The best current result is due to Stevens and de~Zeeuw \cite{SdZ}, who proved that for any set $\P$ of $m$ points in $\mathbb{F}_p^2$ and any set $\L$ of $n$ lines in $\mathbb{F}_p^2$, if  $m^{7/8}<n<m^{8/7}$ and $m^{-2}n^{13}\ll p^{15}$, then the number of incidences between $\P$ and $\L$ satisfies $I(\P, \L)\ll m^{11/15}n^{11/15}.$ Their proof relies mainly on Rudnev's point-plane incidence bound \cite{Rud}. We point out that this theorem also holds in the setting of arbitrary fields $\mathbb{F}$. In this case, the condition $m^{-2}n^{13}\ll p^{15}$ is replaced by $m^{-2}n^{13}\ll \mathtt{char}(\mathbb{F})^{15}$, and is removed if the characteristic of $\mathbb{F}$ is zero. 

In a very recent paper, Rudnev and Wheeler \cite{RJ} obtained point-hyperbola incidence bounds, that is, incidences between points and hyperbolas of the form $(x-a)(y-b)=1$ in $\mathbb{F}_p^2$, which improve earlier results due to Shkredov in \cite{SH21} and Bourgain in \cite{Bourgain}. The main idea in their argument is as follows: for a fixed point $\mathbf{q}\in \P$, let $H_{\mathbf{q}}$ be the set of hyperbolas passing through $\mathbf{q}$. They observed that the number of $k$-rich hyperbolas in $H_{\mathbf{q}}$ can be estimated by using the Stevens--de~Zeeuw point-line incidence bound. This observation has been used by Warren and Wheeler \cite{Audie} to derive an incidence bound between a point set and a set of M\"obius transformations in the plane $\mathbb{F}_p^2$.

The main purpose of this paper is to employ this key idea to study a general incidence problem, namely, incidences between points and irreducible conics in the plane $\mathbb{F}_p^2$. We recall that the three types of irreducible conics include parabolas, ellipses and hyperbolas. We also remark that the results in Sections \ref{Sec:point-conic} and \ref{sec:circ/para} can be extended to arbitrary fields, as long as the relevant characteristic condition is satisfied - this is because each of these results relies on the point line incidence bounds of Stevens and de Zeeuw, which hold over arbitrary fields. The results of Warren and Wheeler also rely on these results, and can therefore also be extended to arbitrary fields.

In Section \ref{sec:applications}, we give various applications of our incidence results. For instance, we prove a pinned algebraic distances result, and an analogue of Beck's theorem for conics in finite fields.

\subsection{Point-conic incidence bounds over prime fields} \label{Sec:point-conic}

Our first result gives an upper bound on the number of incidences between small sets of points and irreducible conics over $\F_p^2$.
\begin{theorem}\label{thm:Conics}
For any set $\mathcal{C}$ of irreducible conics in $\F_p^2$, and any set of points $\mathcal{P} \subseteq \F_p^2$ with $|\mathcal{P}| \ll p^{15/13}$, we have
$$I(\mathcal{P},\mathcal{C}) \ll |\mathcal{P}|^{23/27}|\mathcal{C}|^{23/27} + |\mathcal{P}|^{13/9}|\mathcal{C}|^{12/27} + |\mathcal{C}|.$$
\end{theorem}
To compare with the Cauchy-Schwarz bounds, we note that any conic is determined uniquely by five points, no three of which are collinear, (see \cite[Exercise~69]{Ball}) and by B\'ezout's theorem, any two distinct conics meet in at most four distinct points. So by the K\H{o}v\'{a}ri–S\'{o}s–Tur\'{a}n theorem (see \cite[Theorem~IV.10]{Bol}), we have
\begin{equation}\label{eqn:ConTriv}
    I(\mathcal{P}, \mathcal{C})\ll \min\{|\mathcal{P}||\mathcal{C}|^{4/5} + |\mathcal{C}|,  |\mathcal{P}|^{1/2}|\mathcal{C}| + |\mathcal{P}|\}.
\end{equation}
Theorem~\ref{thm:Conics} improves the trivial bounds of \eqref{eqn:ConTriv}, in the range $|\mathcal{P}|^{19/8} \leq |\mathcal{C}| \leq |\mathcal{P}|^{20/7}$, which encompasses the `balanced Cartesian product' range - that is, when the set of conics $\mathcal{C}$ have coefficients $(a,b,c,d,e)$ coming from a Cartesian product $A \times B \times C \times D \times E$, and the point set $\mathcal{P}$ is also a Cartesian product $F \times G$, and all the sets involved are of the same size $N$. In this case we have $|\mathcal{C}| = N^5 = |\mathcal{P}|^{5/2}$, and our bound improves on \eqref{eqn:ConTriv}.

Next, we give an improvement of Theorem~\ref{thm:Conics} for the particular case when our point set is a Cartesian product.
\begin{theorem}\label{thm:ConicsCP}
Let $\mathcal{C}$ be a set of irreducible conics in $\mathbb{F}_p^2$. Given any sets $A, B\subset \F_p$ with $|A|\leq |B|$ and $|A||\mathcal{C}| \ll p^{2}$, we have 
\[I(A\times B, \mathcal{C})\ll |A|^{\frac{3}{4}}|B|^{\frac{5}{8}}|\mathcal{C}|^{\frac{7}{8}}+|A|^{\frac{1}{2}}|B|^{\frac{3}{4}}|\mathcal{C}|^{\frac{1}{4}} +|\mathcal{C}|.\]
\end{theorem}

In order to compare these results to what is known over the real numbers, the best analogue is given by the Pach-Sharir theorem \cite{PachSharir}, which implies that the number of incidences between a point set $\mathcal{P}$ and an arbitrary set of conics $\mathcal{C}$ satisfies
$$I(\mathcal{P},\mathcal{C}) \ll |\mathcal{P}|^{5/9}|\mathcal{C}|^{8/9} + |\mathcal{P}|+|\mathcal{C}|.$$
{We refer the reader to \cite{sheffer} for a survey of incidence results and their applications over $\R^d$.}

\subsection{Point-circle, point-parabola and point-hyperbola incidence bounds over prime fields} \label{sec:circ/para}
When $\mathcal{C}$ is a set of circles, parabolas or hyperbolas, we have the following improvements. 
\begin{theorem}\label{thm:CircPar}
Let $\mathcal{C}$ be either a set of circles, or of parabolas of the form $y = ax^2 + bx + c$, or of hyperbolas of the form $(x-a)(y-b)=c$ and $\mathcal{P} \subseteq \F_p^2$, with $|\mathcal{P}|\ll p^{15/13}$. If $\mathcal{C}$ is a set of circles suppose that $p\equiv 3 \pmod 4$. Then we have
$$I(P,\mathcal{C}) \ll |\mathcal{P}|^{15/19}|\mathcal{C}|^{15/19} + |\mathcal{P}|^{23/19}|\mathcal{C}|^{4/19} + |\mathcal{C}|.$$
\end{theorem}

We remark that any two parabolas or circles meet in at most two points, and that they are determined uniquely by three non-collinear points, as a consequence of the K\H{o}v\'{a}ri–S\'{o}s–Tur\'{a}n theorem, we have
\begin{equation}\label{eqn:ParCircriv}
    I(\mathcal{P}, \mathcal{C})\ll \min\{|\mathcal{P}||\mathcal{C}|^{2/3} + |\mathcal{C}|,  |\mathcal{P}|^{1/2}|\mathcal{C}| + |\mathcal{P}|\}.
\end{equation}

Theorem~\ref{thm:CircPar} is better than this Cauchy-Schwarz bound in the range $|\mathcal{P}|^{11/8} \ll |\mathcal{C}| \ll |\mathcal{P}|^{12/7}$, again encompassing the balanced Cartesian product range analogous to that described above.

As before, an improved estimate is obtained when $\mathcal{P}$ is a Cartesian product. 
\begin{theorem}\label{thm:CircParCP}
Let $A, B\subset \F_p$, and let $\mathcal{C}$ be either a set of circles, or of parabolas of the form $y = ax^2 + bx + c$, or of hyperbolas of the form $(x-a)(y-b)=c$. If $\mathcal{C}$ is a set of circles, suppose that $p\equiv 3 \pmod 4$. If $|A||\mathcal{C}| \ll p^2$, then we have
$$I(A \times B, \mathcal{C}) \ll |A|^{4/5}|B|^{3/5}|\mathcal{C}|^{4/5} + |A|^{6/5}|B|^{7/5} |\mathcal{C}|^{1/5} + |\mathcal{C}|.$$
\end{theorem}
\subsection{Incidence bounds over arbitrary finite fields for large sets}
Using the same approach, we are able to extend those theorems into arbitrary finite fields for large sets as follows. 
\begin{theorem}\label{arbitaryff}
Let $\mathcal{C}$ be a set of irreducible conics in $\mathbb{F}_q^2$ and $\mathcal{P}$ be a set of points in $\mathbb{F}_q^2$. We have 
\[I(\mathcal{P}, \mathcal{C})\ll \frac{|\mathcal{P}||\mathcal{C}|}{q}+q^{1/5}|\mathcal{P}|^{4/5}|\mathcal{C}|^{4/5}+|\mathcal{C}|.\]
\end{theorem}

\begin{theorem}\label{point-sphere}
Let $\mathcal{P}$ and $\mathcal{S}$ be sets of points and spheres in $\mathbb{F}_q^d$ respectively and assume $q\equiv 3\pmod 4$. We have 
\begin{equation}\label{eqn:PointSphereLS}
    I(\mathcal{P}, \mathcal{S})\ll \frac{|\mathcal{P}||S|}{q}+q^{\frac{d-1}{3}}|\mathcal{P}|^{2/3}|\mathcal{S}|^{2/3}.
    \end{equation}
\end{theorem}

It is worth noting that Theorem \ref{point-sphere} improves earlier results in the literature, namely, it is better than the bound 
\[\frac{|\mathcal{P}||\mathcal{S}|}{q}+q^{\frac{d}{2}}(|\mathcal{P}||\mathcal{S}|)^{1/2}\]
due to Cilleruelo, Iosevich, Lund, Roche-Newton and Rudnev in \cite{CILRR} when $|\mathcal{P}||\mathcal{S}|\le q^{d+2}$, and is better than the bound 
\[\frac{|\mathcal{P}||\mathcal{S}|}{q}+q^{\frac{d-1}{2}}(|\mathcal{P}||\mathcal{S}|)^{1/2},\]
due to Koh, Lee and Pham \cite{KLP}, which holds for small sets $\mathcal{S}$ with $|\mathcal{P}||\mathcal{S}|\le q^{d-1}$.

{In the range $|\mathcal{P}||\mathcal{S}|\ge q^{d+2}$, Theorem \ref{point-sphere} implies the same bound as the above-mentioned result of Cilleruelo et al. in \cite{CILRR}, which is optimal.}

\section{Proof of the point-conic bound (Theorem~\ref{thm:Conics})}

A M\"obius transformation (over $\mathbb F_q$) is a map $f$ of the form
$$f(x) = \frac{ax+b}{cx+d}, \quad ad-bc \neq 0.$$
We will freely swap between the notion of a M\"obius transformation as a map $f$, and the curve given by $y = f(x)$. We require the following result of \cite{Audie} on the number of $k$-rich M\"obius transformations, which are defined as $$T_k := \{ f \in T : |f \cap \mathcal{P}| \geq k \},$$
where $T$ is a set of M\"obius transformations.
\begin{theorem} \label{Mobius}
For any set $T$ of M\"obius transformations, and any set of points $\mathcal{P} \subseteq \F_p^2$ with $|\mathcal{P}| \ll p^{15/13}$, for all $k\geq 3$ we have
$$|T_k| \ll \frac{|\mathcal{P}|^{15/4}}{k^{19/4}} +\frac{|\mathcal{P}|^2}{k^2}.$$
\end{theorem}

\begin{proof}[Proof of Theorem \ref{thm:Conics}]
We shall first aim to bound the number of $k$-rich conics in $\mathcal{C}$, the set of which we denote by $\mathcal{C}_k$. Fix two distinct points $\bm{q_1}$ and $\bm{q_2}$ in $\mathcal{P}$. Define the set $$\mathcal{C}_{\bm{q_1},\bm{q_2},k} := \{C \in \mathcal{C} : \bm{q_1}, \bm{q_2} \in C, |C \cap \mathcal{P}| \geq k\}.$$
$\mathcal{C}_{\bm{q_1},\bm{q_2},k}$ is a set of conics which all pass through the two points $\bm{q_1}$ and $\bm{q_2}$. We apply a projective transformation $\pi$, with the property that $\pi$ maps $\bm{q_1}$ and $\bm{q_2}$ as follows
$$\bm{q_1} \rightarrow [0:1:0], \qquad \bm{q_2} \rightarrow [1:0:0].$$
Note that such a transformation always exists. We now analyse the image $\pi(C)$ for each $C \in \mathcal{C}_{\bm{q_1},\bm{q_2},k}$.

Since $\pi$ is a projective transformation, these images must all remain degree two algebraic curves, that is, conics. Furthermore, we know that the two points at infinity $[0:1:0]$ and $[1:0:0]$ both lie on $\pi(C)$. We claim that this in fact forces $\pi(C)$ to be a M\"obius transformation. Indeed, letting $\gamma$ denote an irreducible conic over $\F_p$, we have
\begin{equation}\label{eqn:MTIFF2RP}
    \gamma \text{ is a M\"obius transformation} \quad \iff \quad \left\{[0:1:0] ,[1:0:0] \right\} \subseteq \gamma.
    \end{equation}

First, we show the implication $\implies$ of \eqref{eqn:MTIFF2RP}. If $\gamma$ is M\"obius, then since it corresponds to a conic, it is given by an equation of the form
$$y = \frac{ax+b}{x+c}$$
which we then rearrange to $$xy + cy - ax - b =0.$$
Projectivising this curve, we have
$$xy + cyz - axz - bz^2 =0.$$
Upon setting $z=0$ to find the points at infinity, we have $xy=0$, and so one of $x$, $y$ must be zero. This yields the two points as in \eqref{eqn:MTIFF2RP}.

We now show the reverse implication. Projectivising an arbitrary conic, of the form
$$ax^2 + by^2 + xy + cy + dx + e = 0,$$
we have 
$$ax^2 + by^2 + xy + cyz + dxz + ez^2=0.$$
Setting $z=0$ to find the points at infinity, we have the equation $ax^2 + by^2 + xy = 0$. We know that the points $(0,1)$ and $(1,0)$ must lie on this curve. The first implies that $b = 0$, and the second that $a =0$. Therefore the original curve $\gamma$ is of the form
$$xy + cy + dx + e = 0.$$
We also must have $dc \neq e$, as otherwise this conic is reducible as $(y+d)(c+x)=0$. Therefore, we have shown that $\gamma$ is the general form of a M\"obius transformation, concluding the proof of \eqref{eqn:MTIFF2RP}.

With \eqref{eqn:MTIFF2RP} at hand, we know that the set $\pi(\mathcal{C}_{\bm{q_1},\bm{q_2},k})$ is a set of M\"obius transformations, and so we can apply Theorem \ref{Mobius}. Therefore, for each distinct pair $(\bm{q_1},\bm{q_2}) \in \mathcal{P}^2$, and for each $k \geq 5$, we have
\begin{equation}\label{krichconic}
    |\mathcal{C}_{\bm{q_1},\bm{q_2},k}| \ll \frac{|\mathcal{P}|^{15/4}}{k^{19/4}} +\frac{|\mathcal{P}|^2}{k^2}.
\end{equation}
Note that the condition $k \geq 3$ has changed to $k \geq 5$, since the two points at infinity on these conics are now being ignored. In order to alter this bound into a bound on $\mathcal{C}_k$, we sum over each distinct pairs of points.
$$|\mathcal{C}_k| \leq {k \choose 2}^{-1}\sum_{\bm{q_1},\bm{q_2} \in \mathcal{P}} |\mathcal{C}_{\bm{q_1},\bm{q_2},k}| \ll \frac{1}{k^2}\sum_{\bm{q_1},\bm{q_2} \in \mathcal{P}} |\mathcal{C}_{\bm{q_1},\bm{q_2},k}|.$$
The binomial factor appears since each $k$-rich conic is being counted at least ${k \choose 2}$ times, once for each pair of distinct points on the conic. We then bound this by ${k \choose 2} \gg k^2$ (this is certainly valid in the range $k \geq 5$). Using \eqref{krichconic}, we have
\begin{equation}\label{krichconic2}
    |\mathcal{C}_k| \ll \frac{|\mathcal{P}|^{23/4}}{k^{27/4}} + \frac{|\mathcal{P}|^4}{k^4}.
\end{equation}
We proceed to use this estimate to obtain the required incidence bound. We use $\mathcal{C}_{=k}$ to denote $\{C\in \mathcal{C}$: $|C\cap \P| = k\}.$

\begin{align*}
    I(\mathcal{P},\mathcal{C}) &= \sum_{k\ge 1} |\mathcal{C}_{=k}|k \\
    & = \sum_{k \leq \Delta} |\mathcal{C}_{=k}|k + \sum_{k > \Delta} |\mathcal{C}_{=k}|k \\
    & \ll \Delta|\mathcal{C}| + \sum_{i\ge 0} \sum_{\substack{C \in \mathcal{C} \\ 2^i \Delta \leq |C \cap \mathcal{P}| < 2^{i+1} \Delta}} (2^{i} \Delta) \\
    & \ll \Delta|\mathcal{C}| + \sum_{i\ge 0}|\mathcal{C}_{2^i\Delta}|(2^{i} \Delta) \\
    & \ll \Delta |\mathcal{C}| + \sum_{i\ge 0}\left(\frac{|\mathcal{P}|^{23/4}}{(2^{i} \Delta)^{27/4}} + \frac{|\mathcal{P}|^4}{(2^{i} \Delta)^4}\right)(2^{i} \Delta) \\
    & \ll \Delta |\mathcal{C}| + \frac{|\mathcal{P}|^{23/4}}{\Delta^{23/4}} + \frac{|P|^4}{\Delta^3}.
\end{align*}

In order to optimise the first two terms, we make the choice
$$\Delta = \max\left\{5, \frac{|\mathcal{P}|^{23/27}}{|\mathcal{C}|^{4/27}}\right\}.$$
This maximum is taken to ensure the application of \eqref{krichconic2} was valid. If the maximum above is $5$, then we must have
$$\frac{|\mathcal{P}|^{23/27}}{|\mathcal{C}|^{4/27}} \leq 5 \implies |\mathcal{P}|^{23} \ll |\mathcal{C}|^{4}.$$
The above bound then becomes
$$I(\mathcal{P},\mathcal{C}) \ll |\mathcal{C}| + |\mathcal{P}|^{23/4} \ll |\mathcal{C}|.$$
If the second term in the maximum is taken, we then have
$$I(\mathcal{P},\mathcal{C}) \ll |\mathcal{P}|^{23/27}|\mathcal{C}|^{23/27} + |\mathcal{P}|^{13/9}|\mathcal{C}|^{12/27}.$$
Summing these two bounds to account for either case then gives the result.
\end{proof}
\section{Proof of the point-conic bound for Cartesian product sets (Theorem \ref{thm:ConicsCP})}
The proof of Theorem \ref{thm:ConicsCP} is similar to that of Theorem \ref{thm:Conics}. First, we recall an incidence result of Stevens and de~Zeeuw. The following version appears in \cite[Theorem~5]{MurPet}. We state the result, more generally, over the two dimensional projective space over $\F_p$, denoted by $\mathbb{P}_2(\F_p)$.

\begin{theorem}\label{thm:ShkPLInc}
Given $A, B \subset \F_p$ with $|A|\leq |B|$, let $\mathcal{P} = \{[a: b: 1]: (a, b)\in A\times B\} \subset \mathbb{P}_2(\F_p)$ and let $\mathcal{L}$ be a set of lines over $\mathbb{P}_2(\F_p)$. Suppose $|A||\mathcal{L}|\ll p^2$. Then
\[
I(\mathcal{P}, \mathcal{L}) \ll |A|^{3/4}|B|^{1/2}|\mathcal{L}|^{3/4} + |\mathcal{L}| + |A||B|.
\]
\end{theorem}

\begin{corollary}\label{Cor:ShkPIPLInc}
Let the sets $\mathcal{P}, \mathcal{L}$ be as in Theorem \ref{thm:ShkPLInc} and let $\pi$ be a projective transformation of $\mathbb{P}_2(\F_p)$. Then
\[
I(\pi(\mathcal{{P}}), \mathcal{L}) \ll |A|^{3/4}|B|^{1/2}|\mathcal{L}|^{3/4} + |\mathcal{L}| + |A||B|.
\]
\end{corollary}

\begin{proof}
First note that $I(\pi(\mathcal{{P}}), \mathcal{L}) = I(\mathcal{{P}}, \pi^{-1}(\mathcal{L}))$. Then the result follows from Lemma~\ref{thm:ShkPLInc} noting that $|\mathcal{L}| = |\pi^{-1}(\mathcal{L})|$.
\end{proof}

\begin{corollary}
For $A, B\subset \F_p$, with $|A|\leq |B|$, let $\mathcal{P} = \{[a: b: 1]: (a, b)\in A\times B\} \subset \mathbb{P}_2(\F_p)$ and let $\mathcal{L}$ denote a set of lines over $\mathbb{P}_2(\F_p)$. Suppose $\pi$ is a projective transformation and for $k\geq 2$, let $\mathcal{L}_k$ denote the set of $k$-rich lines of $\mathcal{L}$ with respect to $\pi(P)$. Suppose $|A||\mathcal{L}|\ll p^2$. Then 
\[
|\mathcal{L}_k|\ll  \frac{|A|^3|B|^2}{k^4} + \frac{|A||B|}{k}.
\]
\end{corollary}

\begin{proof}
The result follows from Corollary~\ref{Cor:ShkPIPLInc}, using the observation that $k|\mathcal{L}_k|\ll I(\pi(\mathcal{P}), \mathcal{L}_k)$.
\end{proof}

By following the arguments of \cite{Audie}, one obtains the following result on the number of $k$-rich M\"obius transformations.
\begin{corollary}\label{cor:kRichTCP}
For $A, B\subset \F_p$, with $|A|\leq |B|$, let $\mathcal{P} = A\times B$ and let $T$ denote a set of M\"obius transformations. Suppose $\pi$ is a projective transformation and for $k\geq 3$, let $T_k$ be a the set of $k$-rich transformations of $T$, with respect to $\pi(P)$. If $|A||T|\ll p^2$, then 
\[
|T_k|\ll \frac{|A|^4|B|^3}{k^5} + \frac{|A|^2|B|^2}{k^2}.
\]
\end{corollary}
\begin{proof}[Proof of Theorem \ref{thm:ConicsCP}]
Corollary~\ref{cor:kRichTCP} may be used, in a similar manner as in the proof of Theorem \ref{thm:Conics}, to bound the number of $k$-rich conics $\mathcal{C}_k$, for $k\geq 5$. This gives
\[
|\mathcal{C}_k|\ll \frac{|A|^6|B|^5}{k^7} + \frac{|A|^4|B|^4}{k^4}.
\]
This can then be converted into the required incidence bound in a similar fashion to the proof of Theorem \ref{thm:Conics}.
\end{proof}

\section{Proof of Theorems \ref{thm:CircPar} and \ref{thm:CircParCP}}
We begin this section by giving the proof of the circles part of Theorem \ref{thm:CircPar}, and then explain the alterations necessary to deal with parabolas and hyperbolas.

\begin{proof}[Proof of Theorem \ref{thm:CircPar}]

Fix a point $\bm{q} \in \mathcal{P}$. We aim to bound the number of $k$-rich circles passing through $\bm{q}$, for $k \geq 3$. After translating the points and circles, we may assume $\bm{q} = (0,0)$ without altering the incidences. We rename the translated sets of points and circles as $\mathcal{P}$ and $\mathcal{C}$ respectively. A circle of the form 
$$(x - c)^2 + (y-d)^2 = r$$
for some $c,d,r$ before this translation, is now of the form 
$$C_{a,b} : \quad (x - a)^2 + (y-b)^2 = a^2 +b^2$$
for some $a,b\in \F_p$; this is due to the fact that $(0,0)$ must lie on the translated circle. Let us take a $k$-rich circle $C$. Then there exist points $(\alpha_1, \beta_1)$,...,$(\alpha_{k-1},\beta_{k-1}) \in \mathcal{P} \setminus \{(0,0)\}$ which all lie on $C$. Each such point can be associated with a line of the form
$$-2\alpha_i X - 2\beta_i Y + \alpha_i^2 + \beta_i^2 = 0.$$

We show that these lines are not defined with multiplicity. Suppose that the two points $(\alpha, \beta)$ and $(\alpha', \beta')$ define the same line. Without loss of generality, we may assume that assume $\beta \neq 0$, since at least one of $\alpha$ or $\beta$ is non-zero. After rearranging, we come to the equation
$$Y = \frac{-\alpha}{\beta}X + \frac{\alpha^2 + \beta^2}{2\beta}.$$
We must therefore have $$\frac{\alpha}{\beta} = \frac{\alpha'}{\beta'}, \quad \frac{\alpha^2 + \beta^2}{2\beta} = \frac{\alpha'^2 + \beta'^2}{2\beta'}$$

setting $\lambda := \frac{\alpha}{\beta} = \frac{\alpha'}{\beta'}$, we must have $\alpha = \lambda \beta$ and $\alpha' = \lambda \beta'$. Substituting this into the second equation then gives
$$\beta(\lambda^2 + 1) = \beta'(\lambda^2 + 1).$$

Recalling our assumption $p\equiv 3\pmod 4$, it follows that $-1$ is a non-square, and so we conclude that $\beta = \beta'$, which implies $\alpha = \alpha'$, and so we are done.

As we have just seen, the point set $\mathcal{P}$ gives rise to a set of lines $\mathcal{L}$, each of the form above, and $|\mathcal{P}| = |\mathcal{L}|$. Furthermore, a circle $C_{a,b}$ as defined above can be associated with the point $(a,b)$. The set of circles $\mathcal{C}$ therefore gives rise to a set of points $\mathcal{Q} \subseteq \F_p^2$. Finally, the $k$-rich circle $C_{a,b}$ above, corresponds to a $(k-1)$-rich point $(a,b) \in \mathcal{Q}$, with respect to the set of lines $\mathcal{L}$. We can see these as the lines corresponding to the points $(\alpha_i,\beta_i)$ all passing the point $(a,b)$. Now, we use the following corollary of the Stevens--de~Zeeuw incidence bound \cite[Theorem~3]{SdZ} to bound the number of such $(k-1)$-rich points. Note that this corollary is the dual of \cite[Corollary~5]{Audie}.
\begin{corollary}\label{cor:RichP}
Let $\mathcal{L}$ be a set of lines over $\F_p^2$, with $|\mathcal{L}|\ll p^{15/13}$, and for $k\geq 2$ let $\mathcal{P}_k$ denote the number of $k$-rich points with respect to $\mathcal{L}$. Then
\[
|\mathcal{P}_k|\ll \frac{|\mathcal{L}|^{11/4}}{k^{15/4}} + \frac{|\mathcal{L}|}{k}.
\]
\end{corollary}

Corollary~\ref{cor:RichP} and the above argument immediately yields the bound
\begin{equation*}
    |\mathcal{C}_{\bm{q}, k}|\ll \frac{|\mathcal{L}|^{11/4}}{k^{15/4}} + \frac{|\mathcal{L}|}{k} = \frac{|\mathcal{P}|^{11/4}}{k^{15/4}} + \frac{|\mathcal{P}|}{k}
\end{equation*}
where $\mathcal{C}_{\bm{q}, k}$ denotes the set of $k$-rich circles, for $k\geq 3$, containing $\bm{q}$. Then summing the contribution over all points of $\mathcal{P}$ and noting that each $k$-rich circle gets overcounted by at least a factor of $k-1$ in this way, we have
\[
|\mathcal{C}_k|\ll \frac{1}{k}\sum_{\bm{q}\in \mathcal{P}}|\mathcal{C}_{\bm{q}, k}| \ll \frac{|\mathcal{P}|^{15/4}}{k^{19/4}} + \frac{|\mathcal{P}|^2}{k^2}.
\]
This can then be used to bound $I(\mathcal{P}, \mathcal{C})$ in a similar manner as in the proofs of Theorem~\ref{thm:Conics} or \cite[Theorem~2]{Audie}.

For parabolas, a similar argument works. After fixing a point $\bm{q}$ and translating so that $\bm{q} = (0,0)$, we find parabolas of the form 
$$y = ax^2 + bx.$$
A $k$-rich parabola then yields a $(k-1)$-rich point $(a,b)$, with respect to the lines given by 
$$\beta = X\alpha^2 + Y\alpha$$
which are again defined with no multiplicity, and the rest of the argument follows similarly to above.

Hyperbolas defined by the equation $(x-a)(y-b)=c$, passing through the point $\mathbf{q}=(q_1, q_2)$, can be written as
\[xy-b(x-q_1)-a(y-q_2)=q_1q_2.\]
By setting $x'=x-q_1$ and $y'=y-q_2$, our hyperbolas are represented by 
\[x'y'+x'(q_2-b)+y'(q_1-a)=0.\]
This equation can be viewed as an incidence between the point $(q_2-b, q_1-a)$ and the line \[x'\cdot X+y'\cdot Y=-x'y'.\]
Thus, we are now in the same situation as before for circles, and the same argument works.  \end{proof}

The proof of Theorem \ref{thm:CircParCP} follows from almost exactly the same argument as that of Theorem~\ref{thm:CircPar}, with the only exception being that the dual form of Theorem~\ref{thm:ShkPLInc} (using point-line duality) replaces Corollary~\ref{cor:RichP} to bound the number of $k$-rich points corresponding to the $k$-rich circles of $\mathcal{C}$. This is based on the observation that the line set $\mathcal{L}$, as defined in the proof of Theorem~\ref{thm:CircPar} (corresponding to the point set $\mathcal{P}$ in the statement of Theorem~\ref{thm:CircParCP}), now has a Cartesian product structure.

\section{Proofs of incidence bounds for large sets (Theorems \ref{arbitaryff} and \ref{point-sphere})}

To prove Theorem \ref{arbitaryff}, we recall the following point-line incidence bound for large sets in \cite{Vinh}. 
\begin{theorem}\label{thm:VinhPL}
Let $\mathcal{P}$ be a set of points and $\mathcal{L}$ be a set of lines in $\mathbb{F}_q^2$. Then 
\[I(\mathcal{P}, \mathcal{L})\le \frac{|\mathcal{P}||\mathcal{L}|}{q}+q^{1/2}\sqrt{|\mathcal{P}||\mathcal{L}|}.\]
\end{theorem}
\begin{corollary} \label{cor:krichlarge}
Let $\mathcal{P}$ be a set of points and $\mathcal{L}_k$ be the set of $k$-rich lines over $\F_q^2$. Suppose that $k>|\mathcal{P}|/q$, then we have 
\[|\mathcal{L}_k|\le \frac{q|\P|}{k^2}.\]
\end{corollary}

Utilising Corollary \ref{cor:krichlarge} and the arguments of \cite[Theorem~2]{Audie}, one obtains the following bound on the number of $k$-rich M\"obius transformations, which is an analogue of Theorem~\ref{Mobius} for large sets.

\begin{corollary}\label{cor:TkLS}
Let $\mathcal{P}$ be a set of points and $T$ a set of M\"obius transformations over $\F_q^2$. Let $T_k$ be the set of $k$-rich transformations in $T$ and suppose that $k> \max\{2, |\mathcal{P}|/q\}$. Then  
\[|T_k|\le \frac{q|\P|^2}{k^3}.\]
\end{corollary}

\begin{proof}[Proof of Theorem \ref{arbitaryff}]

The proof follows the same approach as in the proof of Theorem \ref{thm:Conics}, essentially only replacing Theorem~\ref{Mobius} by Corollary~\ref{cor:TkLS} to bound the number of $k$-rich transformations. In particular, we are able to show that 

\[|\mathcal{C}_k|\ll \frac{|\P|^2}{k^2}\cdot \frac{q|\mathcal{P}|^2}{k^3}=\frac{q|\mathcal{P}|^4}{k^5},\]
whenever $k>\max\{4, |\mathcal{P}|/q\}$. Then, writing
\[
I(\mathcal{P},\mathcal{C}) \le  \sum_{k\leq \frac{|\mathcal{P}|}{q}} |\mathcal{C}_{=k}|k + \sum_{k> \max \{4,\frac{|\mathcal{P}|}{q}\}} |\mathcal{C}_{=k}|k+{\sum_{k\le 4}|\mathcal{C}_{=k}|k},
\]
we may bound the second sum similarly to the proof of Theorem \ref{thm:Conics}, and the first sum trivially, to obtain
\[I(\mathcal{P}, \mathcal{C})\ll \frac{|\mathcal{P}||\mathcal{C}|}{q}+q^{1/5}|\mathcal{P}|^{4/5}|\mathcal{C}|^{4/5}+{|\mathcal{C}|}.\]

\end{proof}

To prove Theorem \ref{point-sphere}, we require the following bound on incidences between large sets of points and hyperplanes due to Vinh~\cite{Vinh}.
\begin{theorem}\label{VinhPHP}
Let $\mathcal{P}$ be a set of points and $\mathcal{H}$ be a set of hyperplanes in $\mathbb{F}_q^d$. The number of incidences between $\mathcal{P}$ and $\mathcal{H}$ satisfies
\[I(\mathcal{P}, \mathcal{H})\le \frac{|\mathcal{P}||\mathcal{H}|}{q}+q^{\frac{d-1}{2}}(|\mathcal{P}||\mathcal{H}|)^{1/2}.\]
\end{theorem}
\begin{proof}[Proof of Theorem \ref{point-sphere}]
The proof uses the same framework as the proof of Theorem \ref{thm:CircPar} and so we skip the overly similar details. We begin by fixing a point $\bm{q}$ aiming to bound the number of $k$-rich spheres passing through it. After a translation, we assume $\bm{q} = \bm{0}$ and so each sphere passing through $\bm{q}$ takes the form
\[
S_{\bm{a}} : \quad (x_1 - a_1)^2 + (x_2 - a_2)^2 + \dots + (x_d - a_d)^2 = a_1^2 + a_2^2 +\dots + a_d^2,
\]
for some $\bm{a} = (a_1, a_2, \dots, a_d)\in \F_q^d$. Let $\bm{q_1}, \dots, \bm{q_{k-1}}$ denote the $k-1$ points on $S_{\bm{a}}$ other than $\bm{0}$ and write $\bm{q_i} = (\alpha_{(i,1)}, \alpha_{(i, 2)},\dots, \alpha_{(i, d)})$. For each $1\leq i\leq k-1$, the point $\bm{q_i}$ can be associated with the hyperplane 
\[
-2\alpha_{(i,1)}X_1 - 2\alpha_{(i,2)} X_2 - \dots - 2\alpha_{(i,d)} X_d + \alpha_{(i,1)}^2 + \alpha_{(i, 2)}^2 +\dots + \alpha_{(i, d)}^2=0.
\]
Arguing similarly as in the proof of Theorem \ref{thm:CircPar}, since $-1$ is a non-square, these hyperplanes are defined without multiplicity.

Now, having established the correspondence between our original sets of points and spheres to new sets of hyperplanes and points respectively, in order to bound the number, $|\mathcal{S}_{\bm{q}, k}|$, of $k$-rich spheres through $\bm{q}$, we require a bound on the number, $|\P_k|$, of $k$-rich points in terms of hyperplanes over $\F_q^d$. To this end, we use Theorem~\ref{VinhPHP}, to obtain
\[
|\P_k| \leq \frac{q^{d-1}|\mathcal{H}|}{k^2}
\]
if $k>|\mathcal{H}|/q$. As in the proof of Theorem~\ref{thm:CircPar}, we use the above estimate to bound the number of $k$-rich spheres, obtaining
\[
|\mathcal{S}_k|\leq \frac{q^{d-1}|\P|^2}{k^3}
\]
if $k>|\mathcal{P}|/q$. This can then be easily converted into the required incidence bound. 
\end{proof}
\section{Applications} \label{sec:applications}
\subsection{Pinned algebraic distances}
Our first application is to the pinned algebraic distance problem. 
\begin{theorem}\label{firstapplication}
Let $f(x, y)$ be one of the following polynomials: $x^2+y^2$ (usual distance function),~ $xy$ (Minkowski distance function) or~ $y+x^2$ (parabola distance function). For $\mathcal{E}\subset \mathbb{F}_p^2$ with $|\mathcal{E}|\ll p^{15/13}$ and $p\equiv 3\pmod 4$, there exists a point $\bm{p}\in \mathcal{E}$ such that $|f(\bm{p}-\mathcal{E})|\gg |\mathcal{E}|^{\frac{8}{15}}$, where 
\[f(\bm{p}-\mathcal{E}):=\{f(\bm{p}-\bm{e})\colon \bm{e}\in \mathcal{E}\}.\]
\end{theorem}

\begin{remark}
When $f(x, y)=x^2+y^2$, Theorem \ref{firstapplication} was first proved by Stevens and de~Zeeuw in \cite{SdZ} by using a point-line incidence bound. The exponent $\frac{8}{15}$ was improved to $\frac{1}{2}+\frac{149}{4214}$ by Iosevich, Koh, Pham, Shen and Vinh \cite{A-P}, then to $\frac{1}{2}+\frac{3}{74}$ by Lund and Petridis \cite{L-P} and to $\frac{1}{2}+\frac{69}{1558}$ by Iosevich, Koh and Pham \cite{IKP}. The best current lower bound is $|\mathcal{E}|^{2/3}$ due to Murphy, Petridis, Pham, Rudnev and Stevens \cite{MPPRP}. We also note that it seems very difficult to extend the methods in \cite{A-P, L-P, IKP, MPPRP} to the Minkowski and parabola distance functions.
\end{remark}
\begin{proof}[Proof of Theorem \ref{firstapplication}]
For $\bm{p}\in \mathcal{E}$, let $\mathcal{C}_{\bm{p}}$ be the set of conics defined by the equation $f(\bm{p}-\bm{x})=t$, with $t\in f(\bm{p}-\mathcal{E})\setminus \{\bm{0}\}$. Let $\mathcal{C}=\bigcup_{\bm{p}\in \mathcal{E}}\mathcal{C}_{\bm{p}}$. We observe that $I(\mathcal{E}, \mathcal{C})\gg |\mathcal{E}|^2$. On the other hand, applying Theorem \ref{thm:CircPar}, we have 
\[I(\mathcal{E}, \mathcal{C})\ll (|\mathcal{E}||\mathcal{C}|)^{15/19}+|\mathcal{E}|^{23/19}|\mathcal{C}|^{4/19}+|\mathcal{C}|.\]
Putting the lower and upper bounds of $I(\mathcal{E}, \mathcal{C})$ together and using the fact that $|\mathcal{C}|\le \sum_{\bm{p}\in \mathcal{E}}|f(\bm{p}-\mathcal{E})|$, the theorem follows.
\end{proof}

We now take advantage of the generality of Theorem \ref{thm:Conics} to study algebraic distances between two sets in $\mathbb{F}_p^3$, where one set lies on a plane and the other set is arbitrary. 

More precisely, let $f\in \F_p[x, y, z]$ and $\mathcal{E}, \mathcal{F}\subset \F_p^3$, with $\mathcal{E}$ lying on the plane $z=0$. As above, the set of $f$-algebraic distances between $\mathcal{E}$ and $\mathcal{F}$ is defined by 
\[f(\mathcal{E}-\mathcal{F}):=\{f(\bm{x}-\bm{y})\colon \bm{x}\in \mathcal{E}, \bm{y}\in \mathcal{F}\}.\]
\begin{theorem}\label{distance-general}
Let $\mathcal{E}, \mathcal{F}\subset \F_p^3$, with $\mathcal{E}$ lying on the plane $z=0$. Suppose $f\in \F_p[x, y, z]$ satisfies the property that, for each $\bm{p}\in \mathcal{F}$, the polynomial $f(\bm{x}-\bm{p})-t$ is of degree two and irreducible for all $t\in f(\mathcal{E}-\bm{p})$. We have 
\[|f(\mathcal{E}-\mathcal{F})|\gg \min \left\lbrace |\mathcal{E}|^{4/23}|\mathcal{F}|^{4/23}, ~\frac{|\mathcal{F}|^{5/4}}{|\mathcal{E}|},~\frac{|\mathcal{E}|^{20/7}}{|\mathcal{F}|}, |\mathcal{E}| \right\rbrace.\]
\end{theorem}
The following is an example of this theorem.
\begin{corollary}
Let $f(x, y, z)=x^2y^2+z^2$. For a set $\mathcal{E}$ on the plane $z=0$ and a set  $\mathcal{F}\subset \mathbb{F}_p^3$ with $|\mathcal{F}|\ge |\mathcal{E}|^{\frac{405}{216}}$ and $p\equiv 3\pmod 4$, we have 
\[|f(\mathcal{E}-\mathcal{F})|\gg \min \left\lbrace |\mathcal{E}|^{4/23}|\mathcal{F}|^{4/23}, ~\frac{|\mathcal{F}|^{5/4}}{|\mathcal{E}|},~\frac{|\mathcal{E}|^{20/7}}{|\mathcal{F}|}, |\mathcal{E}| \right\rbrace.\]
\end{corollary}
\begin{remark}
We note that similar questions for some specific polynomials $f$ have been considered in the literature. For instance, the set of distances between a set on a line and an arbitrary set in $\mathbb{F}_p^2$ was studied by Iosevich, Koh, Pham, Shen and Vinh in \cite{A-P} as the key step in their improvement of Stevens--de Zeeuw's result on the original distance problem in two dimensions. 
\end{remark}
\begin{proof}[Proof of Theorem \ref{distance-general}]
If $|f(\mathcal{E}-\mathcal{F})|>\frac{|\mathcal{E}|^{20/7}}{|\mathcal{F}|}$, then we are done. Thus, assuming otherwise, we have $|f(\mathcal{E}-\mathcal{F})|\cdot |\mathcal{F}|\le |\mathcal{E}|^{20/7}$.

For each $\bm{p}\in \mathcal{F}$, let $\mathcal{C}_{\bm{p}}$ be the set of irreducible conics defined by $f(\bm{x}-\bm{p})=t$ where $t\in f(\mathcal{E}-\bm{p})$. 

Let $\mathcal{C}=\bigcup_{\bm{p}}\mathcal{C}_{\bm{p}}$. We observe that $I(\mathcal{E}, \mathcal{C})\gg|\mathcal{E}||\mathcal{F}|$. On the other hand, applying Theorem~\ref{thm:Conics} implies
\[|\mathcal{E}||\mathcal{F}|\ll |\mathcal{E}|^{\frac{23}{27}}\left( \sum_{\bm{p}\in \mathcal{F}}|f(\mathcal{E}-\bm{p})|  \right)^{23/27}+|\mathcal{E}|^{13/9} \left(\sum_{\bm{p}\in \mathcal{F}}|f(\mathcal{E}-\bm{p})| \right)^{12/27}+\sum_{\bm{p}\in \mathcal{F}}|f(\mathcal{E}-\bm{p})|.\]
Rearranging this inequality gives
\[|f(\mathcal{E}-\mathcal{F})|\gg \min \left\lbrace |\mathcal{E}|^{4/23}|\mathcal{F}|^{4/23}, ~\frac{|\mathcal{F}|^{5/4}}{|\mathcal{E}|}, |\mathcal{E}| \right\rbrace.\]

We note that the assumption $|\mathcal{F}|\ge |\mathcal{E}|^{\frac{405}{216}}$ is needed to ensure that $|\mathcal{E}|^{19/8}\le |\mathcal{C}|\le |\mathcal{E}|^{20/7}$ and so the bound on $I(\mathcal{E}, \mathcal{C})$, given by Theorem~\ref{thm:Conics}, is better than the trivial bound \eqref{eqn:ConTriv}.
\end{proof}

\subsection{Polynomial images}
\begin{theorem}\label{polynomialimage}
Let $f(x, y)$ be either $x^2+y^2$ (usual distance function) or~ $y+x^2$ (parabola distance function). In the former case, assume $p\equiv 3\pmod 4$. For $\mathcal{E}, \mathcal{F}\subset \mathbb{F}_p^2$ with $|\mathcal{E}+\mathcal{F}|\ll p^{15/13}$, we have 
\[|f(\mathcal{E})|\gg \min \left \lbrace \frac{|\mathcal{E}|^{19/15}|\mathcal{F}|^{4/15}}{|\mathcal{E}+\mathcal{F}|}, ~\frac{|\mathcal{E}|^{19/4}|\mathcal{F}|^{15/4}}{|\mathcal{E}+\mathcal{F}|^{23/4}}, |\mathcal{E}|\right\rbrace.\]
\end{theorem}
\begin{proof}
We consider the following equation 
\[f(\bm{x}-\bm{y})=t,\]
where $\bm{x}\in \mathcal{E}+\mathcal{F}$, $\bm{y}\in \mathcal{F}$ and $t\in f(\mathcal{E})$. 

Let $\mathcal{C}$ be the set of curves defined by $f(\bm{x}-\bm{q})=t$ with $\bm{q}\in \mathcal{F}$ and $t\in f(\mathcal{E})$. It is not hard to see that if $f(x, y)=x^2+y^2$ or $f(x, y)=y+x^2$, then the curves in $\mathcal{C}$ are irreducible. Note that
\[|\mathcal{E}||\mathcal{F}|\le I(\mathcal{E}+\mathcal{F}, \mathcal{C}).\]
Since $|\mathcal{C}|=|\mathcal{F}||f(\mathcal{E})|$, using the incidence bound of Theorem~\ref{thm:CircPar}, one has 
\[|\mathcal{E}||\mathcal{F}| \ll |\mathcal{E} + \mathcal{F} |^{15/19}(|\mathcal{F}||f(\mathcal{E})|)^{15/19} + |\mathcal{E} + \mathcal{F}|^{23/19}(|\mathcal{F}||f(\mathcal{E})|)^{4/19} + |\mathcal{F}||f(\mathcal{E})|.\]
Solving this inequality completes the proof. 
\end{proof}

\begin{theorem}\label{polynomialimage2}
Let $f(x, y)=xy$ (Minkowski distance function). For $\mathcal{E}, \mathcal{F}\subset \mathbb{F}_p^2$ with $|\mathcal{E}+\mathcal{F}|\ll p^{15/13}$. Assuming that the two lines $x=0$ and $y=0$ contain at most $|\mathcal{E}|/2$ points from $\mathcal{E}$, we have 
\[|f(\mathcal{E})|\gg \min \left \lbrace \frac{|\mathcal{E}|^{19/15}|\mathcal{F}|^{4/15}}{|\mathcal{E}+\mathcal{F}|}, ~\frac{|\mathcal{E}|^{19/4}|\mathcal{F}|^{15/4}}{|\mathcal{E}+\mathcal{F}|^{23/4}}, |\mathcal{E}|\right\rbrace.\]
\end{theorem}
\begin{proof}
Let $\mathcal{E'}=\mathcal{E}\setminus \{(a, b)\in \mathcal{E}\colon a=0 ~\mbox{or}~ b=0\}$.

By our assumption, we know that $|\mathcal{E}'|\gg |\mathcal{E}|$. We note that $0\not\in f(\mathcal{E}')$, so the curves defined by 
\[(x-a)(y-b)=t,\]
with $(a, b)\in \mathcal{F}$ and $t\in f(\mathcal{E}')$ are irreducible. So, we can use the same argument as in the proof of Theorem \ref{polynomialimage} to conclude the proof of the theorem.
\end{proof}
\begin{remark}
We note that the assumption that there is at most a proportion of points from $\mathcal{E}$ belonging to the two lines $x=0$ and $y=0$ is necessary, for instance, if $\mathcal{E}\subset \{x=0\}\cup \{y=0\}$, then it is clear that $f(\mathcal{E})=\{0\}$.
\end{remark}
\subsection{An improvement of Koh-Sun's result on distances for large sets}
For $\mathcal{E}, \mathcal{F}\subset \mathbb{F}_q^d$, we denote the set of distances between $\mathcal{E}$ and $\mathcal{F}$ by the set $\Delta(\mathcal{E}, \mathcal{F})$. 

We recall results of Koh and Sun~\cite{koh}, which remove the logarithmic factor in a result due to Dietmann \cite{Di12}. In \cite[Theorems~3.3]{koh}, the authors prove that if $d\geq 3$ is odd, then
\begin{equation}\label{main1}
|\Delta(\mathcal{E},\mathcal{F})|\geq \left\{ \begin{array}{ll} \min \left\{ \frac{q}{2}, \frac{|\mathcal{E}||\mathcal{F}|}{8q^{d-1}}\right\}
\quad &\mbox{if} ~~ 1\leq |\mathcal{E}|< q^{\frac{d-1}{2}}\\
\min \left\{ \frac{q}{2}, \frac{|F|}{8q^{\frac{d-1}{2}}}\right\}\quad &\mbox{if} ~~ q^{\frac{d-1}{2}}
\leq |E|< q^{\frac{d+1}{2}}\\
\min \left\{ \frac{q}{2}, \frac{|\mathcal{E}||\mathcal{F}|}{2q^{d}}\right\}\quad &\mbox{if} ~~q^{\frac{d+1}{2}} \leq
|\mathcal{E}|\leq q^{d} \end{array}\right. .
\end{equation}
For even $d\geq 2$, under the assumption $|\mathcal{E}||\mathcal{F}|\geq 16 q^d,$ by \cite[Theorems~3.5]{koh}, one has
\begin{equation} \label{main2}
 |\Delta(\mathcal{E},\mathcal{F})|\geq \left\{
\begin{array}
{ll} \frac{q}{144} \quad&\mbox{for}~~1\le |\mathcal{E}|< q^{\frac{d-1}{2}}\\
\frac{1}{144} \min\left\{ q, \frac{|\mathcal{F}|}{ 2q^{\frac{d-1}{2}}}\right\} \quad&\mbox{for}~~  q^{\frac{d-1}{2}}\leq |\mathcal{E}|< q^{\frac{d+1}{2}}\\
\frac{1}{144}     \min\left\{ q, \frac{2|\mathcal{E}||\mathcal{F}|}{q^d}\right\} \quad&\mbox{for}~~  q^{\frac{d+1}{2}}\leq
|\mathcal{E}|\leq q^{d}
\end{array}\right. .
\end{equation}
We mention that, in comparison to estimate \eqref{main1} (for odd $d$), the additional condition $|\mathcal{E}||\mathcal{F}|\ge 16q^d$ for estimate \eqref{main2} (for even $d$) is necessary in Koh and Sun's proof. This is due to the fact that the Fourier decay of the sphere of zero radius in even dimensions is much worse than in odd dimensions. We also note that in the range $q^{\frac{d+1}{2}}\le |\mathcal{E}|\le q^d$, the lower bound $\gg \min\left\{ q, \frac{|\mathcal{E}||\mathcal{F}|}{q^d}\right\}$ was obtained by Shparlinski \cite{Shpar} without the condition $|\mathcal{E}||\mathcal{F}|\gg q^d$.

As a direct consequence of Theorem \ref{point-sphere}, we are able to remove the condition $|\mathcal{E}||\mathcal{F}|\gg q^d$ for the range $q^{\frac{d-1}{2}}\le |\mathcal{E}|\le q^{\frac{d+1}{2}}$. 

\begin{theorem}
Let $\mathcal{E}, \mathcal{F}$ be sets in $\mathbb{F}_q^d$. Assume that $|\mathcal{E}|\sim |\mathcal{F}|\le q^{\frac{d+1}{2}}$, then we have 
\[|\Delta(\mathcal{E}, \mathcal{F})|\gg \min \left\lbrace q, \frac{|\mathcal{E}|^{1/2}|\mathcal{F}|^{1/2}}{q^{\frac{d-1}{2}}} \right\rbrace.\]
\end{theorem}
\begin{proof}
For $\bm{p}\in \mathcal{F}$, let $\mathcal{C}_{\bm{p}}$ be the set of spheres centered at $\bm{p}$ of radius in $\Delta(\bm{p}, \mathcal{E})$, and let $\mathcal{C}=\bigcup_{\bm{p}\in \mathcal{F}}\mathcal{C}_{\bm{p}}$. We have $|\mathcal{C}|=\sum_{\bm{p}\in \mathcal{F}}|\Delta(\bm{p}, \mathcal{E})|$, and $I(\mathcal{E}, \mathcal{C})\gg |\mathcal{E}||\mathcal{F}|$. Therefore, applying Theorem~\ref{point-sphere} gives 
\[|\mathcal{E}||\mathcal{F}|\ll \frac{|\mathcal{E}|\sum_{\bm{p}\in \mathcal{F}}|\Delta(\bm{p}, \mathcal{E})|}{q}+q^{\frac{d-1}{3}}|\mathcal{E}|^{2/3}\left(\sum_{\bm{p}\in \mathcal{F}}|\Delta(\bm{p}, \mathcal{E})|\right)^{2/3}.\]
Solving this inequality, there exists $\bm{p}\in \mathcal{F}$ such that \[|\Delta(\bm{p}, \mathcal{E})|\gg  \min \left\lbrace q, \frac{|\mathcal{E}|^{1/2}|\mathcal{F}|^{1/2}}{q^{\frac{d-1}{2}}} \right\rbrace.\]
\end{proof}

\subsection{Conical Beck's Theorem}
A well-known result of Beck~\cite{Beck} states that for a finite point set $\mathcal{P}\subset \R^2$, either a positive proportion of $\P$ is collinear or there exist $\Omega(|\P|^2)$ distinct lines supported on $2$-tuples of (distinct) points of $\P$. See \cite[Corollary 14]{SdZ} for a quantitatively weaker analogue of this result, which holds over arbitrary fields. 

We proceed to prove a finite field analogue of Beck's theorem, replacing the notion of lines by irreducible conics. See also \cite[Corollary 2]{Audie} for a similar result involving M\"obius transformations. In the following, an irreducible conic is said to be defined by a point set $\P$ if it passes through at least five points of $\P$.

\begin{theorem}\label{thm:ConBeck}
Let $\P \subseteq \F_p^2$ be a set of points with $|\P| \ll p^{15/13}$, with no positive proportion of $\P$ collinear. Then either there exists a conic $C$ such that $|C \cap \P| \gg |\P|$, or $\P$ defines at least $|\P|^{20/7}$ irreducible conics.
\end{theorem}
We note that this theorem implies that any point set $\P \subseteq \F_p^2$ with $p^{1+\varepsilon}<|\P| \ll p^{15/13}$, for some $\varepsilon>0$, must define at least $|\P|^{20/7}$ irreducible conics, since no conic or line can contain a positive proportion of this point set. We further point out that proof of this theorem relies on a bound on the number of $k$-rich conics obtained as part of the proof of Theorem~\ref{thm:Conics}. Moreover, note that a collinearity restriction is necessary in this theorem; an irreducible conic is only defined by five points which lie in general position. Through essentially the same scheme, one may replace this bound by the one obtained in the proof of Theorem~\ref{thm:CircPar} to give an improved result, concerning circles, parabolas and hyperbolas.
\begin{proof}[Proof of Theorem~\ref{thm:ConBeck}]
We begin by claiming that a positive proportion of the 5-tuples defined by $\P$ are in general position. Indeed, if we let $L$ be the maximum number of collinear points in $\P$, we have
$$\#\text{5-tuples in GP} = |\P|(|\P|-1)(|\P| - L)(|\P| - 3L)(|\P| - 6L) = |\P|^5 - O(L|\P|^4) \gg |\P|^5,$$
where the last inequality follows from $L = o(|\P|)$. From this point we focus only on 5-tuples of points from $\P$ which are in general position.

With a slight abuse of notation, we write $\mathcal{C}_k$ to denote the set of irreducible conics over $\F_p^2$ containing at least $k$ and at most $2k-1$ points of $\P$. The number of $5$-tuples of points of $\P$ contained in elements of $\mathcal{C}_k$ is at most $O(k^5|\mathcal{C}_k|)$ which, by \eqref{krichconic2}, is bounded by $O(|\P|^{23/4}k^{-7/4} + |\P|^4k)$. 

Let $I = \{k\geq 5: \lambda|\P|^{3/7} \leq k \leq \lambda^{-7/4}|\P|\}$ for some constant $\lambda>0$ to be determined. Then, by the above estimate, the total number of $5$-tuples of $\P$ contained in the conics $\cup_{k\in I}\mathcal{C}_k$ is at most $O(\lambda^{-7/4}|\P|^5)$. Taking $\lambda$ to be sufficiently large, we may assume these account for less than, say, half of the total number of $5$-tuples of $\P$. Moreover, we may assume there exists no irreducible conic containing more than $\lambda^{-7/4}|\P|$ points of $\P$, since otherwise there is nothing to prove. 

We conclude that a positive proportion of the $5$-tuples of $\P$ lie on conics belonging to the set 
\[\mathcal{C}:= \bigcup_{k<\lambda|\P|^{3/7}} \mathcal{C}_k.\]
Consequently, we have
\[
|\P|^5 \approx \sum_{C\in \mathcal{C}}|C\cap \P|^5 \ll |\mathcal{C}||\P|^{15/7},
\]
which gives the second possibility claimed by the theorem, concluding the proof.
\end{proof}

\section*{Acknowledgements}
Thang Pham would like to thank to the VIASM for the hospitality and for the excellent working condition. Audie Warren was supported by Austrian Science Fund FWF grant P-34180. We thank Oliver Roche-Newton for helpful discussions, and two reviewers for many valuable comments.

\end{document}